\documentclass[a4paper]{article}
\usepackage{amsmath,amsthm,amssymb,graphicx,multirow,booktabs}
\usepackage[left=3cm, right=2cm, top=1.5cm, bottom=2.2cm]{geometry}
\usepackage{algorithm}
\usepackage[noend]{algpseudocode}
\usepackage{subfig}
\usepackage{longtable,color,float}
\usepackage{multirow}
\usepackage[unicode,hidelinks]{hyperref}
\providecommand{\keywords}[1]{\textit{Keywords:} #1}
\usepackage[numbers,sort&compress]{natbib}
\usepackage{siunitx}
\begin{document}
\title{An open-source solver for finding global solutions to constrained derivative-free optimization problems}

\author{Gannavarapu Chandramouli\\
Industrial Engineering and Operations Research, IIT Bombay Mumbai 400076, India\\
\and
Vishnu Narayanan\\
Industrial Engineering and Operations Research, IIT Bombay Mumbai 400076, India\\
}
\maketitle

\abstract{%
	In this work, we propose a heuristic based open source solver for finding global solution to constrained derivative-free optimization (DFO) problems.
Our solver named Global optimization using Surrogates for Derivative-free Optimization (GSDO) relies on surrogate approximation to the original problem.
In the proposed algorithm, an initial feasible point is first generated.
This point is subsequently used to generate well spaced feasible points for formulating better radial basis function based surrogate approximations to original objective and constraint functions.
Finally, these surrogates are used to solve the derivative-free global optimization problems.
The proposed solver is capable of handling quantifiable and nonquantifiable as well as relaxable and unrelaxable constraints.
We compared the performance of proposed solver with state of the art solvers like Nonlinear Optimization by Mesh Adaptive Direct Search (NOMAD), differential evolution (DE) and Simplicial Homology Global Optimization (SHGO) on standard test problems.
The numerical results clearly demonstrate that the performance of our method is competitive with respect to other solvers.
}%


\maketitle

\section{Introduction}
Many real world applications are modeled as complex simulation problems which involve optimization with expensive objective and constraint functions.
For such problems, in addition to expensive function evaluations, the derivative information of the objective and constraint functions is either unavailable or too expensive to compute.
Such optimization problems fall under the purview of constrained DFO (\cite{conn2009introduction, larson_menickelly_wild_2019, AuHa2017a}).
The use of standard numerical gradient based methods is limited for these problems as they involve a lot of expensive function evaluations, and in some cases, may also be inaccurate. 
Further, automatic differentiation is not applicable due to inherent complexity in the evaluation of functions.
Another major challenge in this domain is the requirement of finding a global solution to the given problem which involves significantly more computational effort in comparison to its local counterpart. 
Evolutionary algorithms (\cite {storn1997differential,endres2018simplicial}) can be used to solve similar problems, albeit, at the cost of large number of function evaluations.

The mathematical formulation for this optimization problem is given as:
\begin{equation}
	\tag{P1}
	\label{prob:originalProb}
	\begin{aligned}
		&\text{minimize } f(x)\\
		&\text{s.t. } g_j(x) \geq 0, \qquad j=1,\ldots, m\\
		&x_i \in [l_i,u_i], \qquad \qquad i=1,\ldots, d
	\end{aligned}
\end{equation}
where $d$ and $m$ denote the dimension and number of constraints of the problem, and $x \in \mathbb{R}^d$.
Also, $l_i$ and $u_i$ denote the finite lower and upper bounds for the $i^{th}$ variable.
Further, $f:\mathbb{R}^n \rightarrow \mathbb{R}$ and $g_j:\mathbb{R}^n \rightarrow \mathbb{R}$ denote the computationally expensive functions for the objective and the $j^{th}$ constraint.
The feasible region defined solely by the bounds constraints is denoted as $\Omega_B \subset \mathbb{R}^d$.
On the other hand, we denote the feasible region defined by the expensive constraints along with bound constraints as $\Omega \subset \mathbb{R}^d$ and assume that $\text{int}(\Omega)$ is not empty.
In this work, we consider the following types of constraints as defined in the constraint taxonomy by \cite{digabel2015taxonomy}:
\begin{enumerate}
	\item Quantifiable, Relaxable, Simulation and Known (QRSK).
	\item Nonquantifiable, Relaxable, Simulation and Known (NRSK).
	\item Quantifiable, Unrelaxable, Simulation and Known (QUSK).
	\item Nonquantifiable, Unrelaxable, Simulation and Known (NUSK).
	\item Nonquantifiable, Unrelaxable, Simulation and Hidden (NUSH).
\end{enumerate}

The different attributes used in above constraints, as outlined by \cite{digabel2015taxonomy}, are briefly described as follows: (1) \emph{quantifiable}: extent of feasibility and/or violation is quantifiable for the given constraint, (2) \emph{relaxable}: output generated by the simulation is meaningful and trustworthy even if the constraint is violated, (3) \emph{simulation}: feasibility of this constraint is verified only by running a simulation and (4) \emph{known}: knowledge of the constraint is explicitly provided in the problem formulation.

In DFO, pattern search (\cite{conn2009introduction, AuHa2017a}) is a prominent approach which uses a predefined geometric pattern to search for a better point around the current solution.
In literature, several pattern search methods have been proposed such as line search (\cite{fasano2014linesearch}), direct search (\cite{kolda2003optimization}) and mesh adaptive direct search (\cite{AuDe2006}).
The other major approach in DFO includes surrogate approximation model (\cite{forrester2008engineering,talgorn2015statistical}).
Surrogate models do not use any predefined pattern but rather try to approximate the actual expensive function. 
The approximated function is cheaper to evaluate than the original one, and its derivative information can be easily obtained.
Thus, standard gradient based approaches can be readily deployed to solve the optimization problem comprising only of the surrogate approximates.
The prominent surrogate models include radial basis functions (\cite{gutmann2001radial,wild2013global, muller2017gosac}), kriging (\cite{kleijnen2009kriging,jones1998efficient}), polynomial regression models (\cite{myers2016response}) and multivariate adaptive regression splines (\cite{friedman1991multivariate}). 
To obtain the global solution to a DFO problem, synchronization between the exploration and exploitation steps within the feasible region is needed, and related investigations have been carried out in several works by \cite{muller2013so,regis2007stochastic,gutmann2001radial} and \cite{jones1998efficient}. 
The exploration step involves search for new points in the unexplored regions of feasible space.
The exploitation step involves search for best objective function value within the feasible space using the current surrogate model.
For more information, an exhaustive review of different DFO solvers is presented by \cite{rios2013derivative}.

Finding the global solution to a constrained DFO problem is a challenging task, especially in the presence of a tight budget on the number of function evaluations.
Solvers such as RBFOpt (\cite{costa2018rbfopt}), EGO (\cite{jones1998efficient}), DIRECT (\cite{jones1993lipschitzian}), DYCORS (\cite{regis2013combining}), AQUARS (\cite{regis2013quasi}), SHEBO (\cite{muller2019surrogate}) and SNOBFIT (\cite{huyer2008snobfit} have been proposed, but they are primarily limited to handle the bound constraints.
For constrained optimization with infeasible starting points, algorithms like COBRA (\cite{regis2014constrained}) and KCGO (\cite{li2017kriging}) have been reported, but their test problems were mainly limited to QRSK type constraints, and to the best of our knowledge, their source codes are not openly available.
Evidently, there seems to be a scarcity of free and open source software (FOSS) libraries which are designed to address global DFO problems with general constraints.
Thus, the main goal of our work is to develop a new solver to handle this challenging problem within limited budget.
Accordingly, we propose an open source solver called Global optimization using Surrogates for Derivative-free Optimization problems (GSDO).

The paper is organised as follows.
In Section \ref{rbf_theory}, we provide a brief overview of radial basis function (RBF).
We describe the details of our proposed approach in Section \ref{propApp}. 
We report the results of our computational experiments in Section \ref{numerical}, followed by software details in Section \ref{software} and conclusions in Section \ref{conclusion}.

\section{Radial Basis Function}\label{rbf_theory}
Given a function $f:\mathbb{R}^d \rightarrow \mathbb{R}$, and a set of $n$ distinct points $\{x^1,\ldots,x^n \} \subseteq \mathbb{R}^d$ with function values $\{f(x^1),\ldots,f(x^n)\}$, we define RBF interpolant (\cite{muller2017gosac}) at $x$ as
\begin{equation}
	\label{rbf_interpolant}
	s_{\text{rbf}}(x) = \sum\limits^n_{i=1}\gamma_i \phi(||x-x^i||_2) + p(x),
\end{equation}
where $\phi:\mathbb{R}_+ \rightarrow \mathbb{R}$, $\gamma_1,\ldots,\gamma_n \in \mathbb{R}$ and $p$ is a polynomial of degree $t$, whose value is dependent upon the underlying form of the function $\phi$.

In this work, cubic radial basis function was used for the surrogate approximation.
Our choice of this surrogate function is primarily motivated by its extensive use in the domain of DFO, as documented by \cite{wild2013global,bjorkman2000global} and \cite{regis2011stochastic}.
For cubic RBF, with $t$ set to $1$, we have the following form of RBF interpolant:
\begin{equation}
	\label{rbf_cubic}
	s_{\text{rbf}}(x) = \sum\limits^n_{i=1}\gamma_i (||x-x^i||_2^3) + \lambda^T \binom{x}{1},
\end{equation}
where $\lambda \in \mathbb{R}^{d+1}$.
The RBF parameters $\gamma$ and $\lambda$ are obtained by solving the following system of linear equations:
\begin{equation}
	\label{eq:rbf_linear_system}
	\begin{bmatrix}
		\Phi & P \\
		P^T & 0_{(d+1)\times(d+1)}\\
	\end{bmatrix}
	\begin{bmatrix}
		\gamma\\
		\lambda\\ 
	\end{bmatrix}
	=
	\begin{bmatrix}
	F\\
	0_{d+1}\\
	\end{bmatrix}
	,
\end{equation}
where $\Phi_{ij} = \phi(||x^i-x^j||)$ for $i,j=1,\ldots,n$, $F=[f(x^1),\ldots,f(x^n)]^T$ and $P$ is a matrix defined as:
\begin{equation}
	\label{eq:rank_rbf}
	P=
	\begin{bmatrix}
		(x^1)^T &1\\
		(x^2)^T &1\\
		\vdots &\vdots\\
		(x^n)^T &1\\
	\end{bmatrix}
\end{equation}
As per \cite{powell1992theory} and \cite{muller2017gosac}, we have the linear system (\ref{eq:rbf_linear_system}) solvable if and only if the rank of $P$ is $d+1$.
Thus, the number of distinct points ($n$) required for RBF model must be greater than or equal to $d+1$.
The above RBF model is used to approximate the objective and constraint functions in the proposed approach, which is discussed next.

\section{Proposed Approach}\label{propApp}
\begin{table}
\caption{List of Symbols}
	\begin{tabular}{|c|p{15cm}|}
	\hline
	$\mathcal{G}$ & Constraint vector $\mathcal{G}:\mathbb{R}^d \rightarrow \mathbb{R}^m$ such that $\mathcal{G}(x) = (g_1(x),\ldots,g_m(x))^T$.\\
	$g^s$ & Surrogate model function for constraint $g$.\\
	$m'$ & Number of quantifiable constraints (QRSK and QUSK).\\
	$\mathcal{G}^s$ & Surrogate constraint vector $\mathcal{G}^s:\mathbb{R}^d \rightarrow \mathbb{R}^{m'}$ such that $\mathcal{G}^s(x) = (g_1^s(x),\ldots,g_{m'}^s(x))^T$.\\
	$g_c$ & Classification constraint for handling all nonquantifiable constraints.\\
	$\mathcal{X}$ & Set of points at which objective function and constraints have been evaluated.\\
	$\mathcal{X}_F$ & Set of feasible points in $\mathcal{X}$.\\
	$\mathcal{X}_I$ & Set of points in $\mathcal{X}$ violating only QRSK constraints but not any of the NRSK, QUSK, NUSK or NUSH constraints.\\
	$\mathcal{X}_S$ & Set of points in $\mathcal{X}$ violating NRSK constraints but not any of the QUSK, NUSK or NUSH constraints.\\
	$\mathcal{X}_U$ & Set of points in $\mathcal{X}$ violating only one QUSK constraint but not any of the NUSK or NUSH constraints.\\
	$\mathcal{X}_H$ & Set of points in $\mathcal{X}$ violating at least one NUSH or NUSK constraint, or multiple QUSK constraints.\\
	$T_{LH}$ & Maximum number of points to be generated using Latin Hypercube (LH) search.\\
	$T_{max}$ & Maximum budget for the problem.\\
	$B(x,r)$  & A ball of radius $r$ with center $x$.\\
	$\eta_{max}$ & Maximum number of feasible points to be generated in Stage-2.\\
	$K_{max}$ & Maximum number of iterations for Stage-2.\\
	$\Omega$ & Feasible search space.\\
	$\Omega_B$ & Feasible space defined by bound constraints.\\
	$\Delta$ & Maximum distance of a point in $\Omega_B \backslash \mathcal{X}_F$ with respect to the points in $\mathcal{X}_F$. \\
	$\Delta_{min}$ & Minimum allowed value of $\Delta$.\\
	$\delta_{lu}$ & Minimum distance between bounds on the variable. So, $\delta_{lu}=\underset{i=1,\ldots,d}{\min} (|u_i-l_i|)$ and $\delta_{lu} > 0$ as $\text{int}(\Omega) \neq \emptyset$.\\
	$|\mathcal{X}|$ & Number of points in the set $\mathcal{X}$.\\
	$\text{M}(\mathcal{X})$ & Matrix of points in the set $\mathcal{X}$ arranged in column-wise order.\\
	\hline
\end{tabular}
\label{table_symb}
\end{table}

The list of symbols and notations used in the proposed method are outlined in Table \ref{table_symb}.
We first introduce few auxiliary Algorithms (\ref{algo:filter}-\ref{algo:random_pt}) and then discuss the primary Algorithms (\ref{algo:feasible}-\ref{algo:global}) which utilize them.
These auxiliary algorithms focus on how to classify a point based on its constraint violation, how to create surrogate models for the quantifiable constraints, how to address the nonquantifiable constraints, and lastly, how to handle situations where the surrogate models are infeasible.

\begin{algorithm}[h]
	\caption{Filter Point $x \in \mathbb{R}^d$}
	\label{algo:filter}
	\begin{algorithmic}[1]
		\State \textbf{If} {$x \notin \mathcal{X}$} \textbf{then} add $x$ to set $\mathcal{X}$.
		\State \textbf{If} {$x$ violates any NUSH, NUSK or at least two QUSK constraints} \textbf{then} add $x$ to set $\mathcal{X}_H$.
		\State \textbf{else if} {$x$ point violates only one QUSK constraint} \textbf{then} add $x$ to set $\mathcal{X}_U$.
		\State \textbf{else if} {$x$ point violates any NRSK constraint} \textbf{then} add $x$ to set $\mathcal{X}_S$.
		\State \textbf{else if} {$x$ point violates only QRSK constraints} \textbf{then} add $x$ to set $\mathcal{X}_I$.
		\State \textbf{else} add $x$ to set $\mathcal{X}_F$.
	\end{algorithmic}
\end{algorithm}

\begin{algorithm}[h]
	\caption{Check the rank of input points for RBF}
	\label{algo:rbf_rank}
	\begin{algorithmic}[1]
		\State \textbf{If} {$|\mathcal{X}| = 0$} \textbf{then} return FALSE. \textbf{Else}, set $k=|\mathcal{X}_F| + |\mathcal{X}_I| + |\mathcal{X}_S|$.
		\State Create matrix $P \in \mathbb{R}^{d\times k}$ such that $P = [M(\mathcal{X}_F) \quad M(\mathcal{X}_I) \quad M(\mathcal{X}_S)]$.
		\For{each QRSK or QUSK constraints}
			\State \textbf{If} {constraint is QUSK} \textbf{then} include points (which are feasible to current constraint) from $\mathcal{X}_U$ into $P$ i.e. $P = [P \quad M(\mathcal{X}_U)]$.
			\State Create a vector $e\in \mathbb{R}^k$ of ones where $k$ equals the number of columns in $P$.
			\State Find rank $r$ of $[P^T \quad e]$.
			\State \textbf{If} {$r < d+1$} \textbf{then} return FALSE.
			\State Restore $P$ such that $P = [M(\mathcal{X}_F) \quad M(\mathcal{X}_I) \quad M(\mathcal{X}_S)]$.
		\EndFor
		\State \textbf{If} rank condition is satisfied for all QRSK and QUSK constraints \textbf{then} return TRUE.
	\end{algorithmic}
\end{algorithm}

\begin{algorithm}[h]
	\caption {Fitting Surrogate Model for a Given Constraint $g$}
	\label{algo:surrogate_fit}
	\begin{algorithmic}[1]
		\State Collect all points from $\mathcal{X}_F$, $\mathcal{X}_I$ and $\mathcal{X}_S$ into set $\mathcal{M}$.
		\State \textbf{If} {constraint $g$ is of QRSK type} \textbf{then} fit RBF model over points in $\mathcal{M}$.
		\If {constraint $g$ is of QUSK type}
			\State Filter points from $\mathcal{X}_U$ which are feasible for the constraint $g$.
			\State Add these filtered points to set $\mathcal{M}$.
			\State Fit RBF model over points in $\mathcal{M}$.
		\EndIf
	\end{algorithmic}
\end{algorithm}

\begin{algorithm}[h]
	\caption{Computing Classification Constraint ($g_{c}(x)$)}
	\label{algo:knn}
	\begin{algorithmic}[1]
		\State \textbf{Input}: $x \in \mathbb{R}^d$.
		\State \textbf{Parameters}: Real numbers $c_1, c_2, c_3, c_4$ such that $c_1 > 0$ and $c_4 < c_3 < c_2 < 0$.
		\If {($\mathcal{X}_F > 2$ or $\mathcal{X}_I > 2$) and ($\mathcal{X}_H > 2$ or $\mathcal{X}_U > 2$ or $\mathcal{X}_S > 2$)}
			\State Use K-Nearest Neighbor Algorithm to classify $x$ into one among the sets $\mathcal{X}_F$, $\mathcal{X}_I$, $\mathcal{X}_S$, $\mathcal{X}_U$ and $\mathcal{X}_H$.
			\State \textbf{If} {$x$ is classified into $\mathcal{X}_F$ or $\mathcal{X}_I$} \textbf{then} return $c_1$.
			\State \textbf{else if} {$x$ is classified into $\mathcal{X}_S$} \textbf{then} return $c_2$.
			\State \textbf{else if} {$x$ is classified into $\mathcal{X}_U$} \textbf{then} return $c_3$.
			\State \textbf{else} \textbf{then} return $c_4$.
		\Else
			\State Return $c_1$.
		\EndIf
	\end{algorithmic}
\end{algorithm}

\begin{algorithm}[h]
	\caption{Create Random Point}
	\label{algo:random_pt}
	\begin{algorithmic}[1]
		\State \textbf{Parameters}: $\delta_r > 0$ and $\delta_d > 0$
		\State Pick the latest evaluated feasible point, $x$, from the set $\mathcal{X}_F$.
		\State Set radius $r = \text{min}\{\delta_{lu}/\delta_r, \delta_d\sqrt{d}\}$.
		\State Create a random point $y \in B(x,r)$ such that $y \in \Omega_B$ and $y \notin \mathcal{X}$.
		\State Return $y$.
	\end{algorithmic}
\end{algorithm}

\subsection{Classification of an evaluated point}
Given a point $x\in \mathbb{R}^d$ and its corresponding constraint values $\mathcal{G}(x)$, Algorithm (\ref{algo:filter}) is applied to classify it into mutually exclusive sets $\mathcal{X}_F$, $\mathcal{X}_I$, $\mathcal{X}_S$, $\mathcal{X}_U$ and $\mathcal{X}_H$.
These sets are defined in Table \ref{table_symb}.
Note that $\mathcal{X} = \mathcal{X}_F \cup \mathcal{X}_I \cup \mathcal{X}_S \cup \mathcal{X}_U \cup \mathcal{X}_H$.
This classification is helpful in building the surrogate models over quantifiable constraints, and also forming a new constraint to explicitly handle the nonquantifiable constraints.
We now define an approach to build the surrogate models for quantifiable constraints using points from the sets $\mathcal{X}_F$, $\mathcal{X}_I$, $\mathcal{X}_S$ and $\mathcal{X}_U$.

\subsection{Handling quantifiable constraints using surrogate model}
Algorithm (\ref{algo:rbf_rank}) is used to ensure that for each quantifiable constraint, the necessary conditions for solvability of Equation (\ref{eq:rbf_linear_system}) are met.
Essentially, it constructs the matrix $P$ for each such constraint as described in the Equation (\ref{eq:rank_rbf}) and checks whether the corresponding rank is equal to $d+1$.
Next, the surrogate models are constructed for QRSK and QUSK constraints using Algorithm (\ref{algo:surrogate_fit}).
Only the quantifiable constraints are handled through the surrogate models.
For nonquantifiable constraints, we applied a different approach as outlined below.

\subsection{Handling nonquantifiable constraints}
We define a new constraint, denoted by $g_c$, to handle all the nonquantifiable constraints at once.
This essentially relies on K-Nearest Neighbor classification approach (\cite{murphy2012machine}).
If a point is classified into the sets $\mathcal{X}_F$ or $\mathcal{X}_I$, then this constraint returns a positive value.
On the other hand, if it is classified into the sets $\mathcal{X}_S$, $\mathcal{X}_U$ or $\mathcal{X}_H$, then it returns a negative value whose magnitude is dependent on the type of classification.
Essentially, this constraint enforces an optimization algorithm to generate points away from the boundary of hidden constraints.
This approach is described in Algorithm (\ref{algo:knn}).

\subsection{Creation of random point}
In spite of having good surrogate models and methods to handle non-quantifiable constraints, it is possible that their corresponding optimization problems can either become infeasible or fail to create a unique point different from the points in the set $\mathcal{X}$.
In order to handle such situations, we use the Algorithm (\ref{algo:random_pt}) to create a random point within feasible region $\Omega_B$.
Also, this approach can be helpful in creating better surrogate models and classification constraints by creating a new point in the unexplored regions of $\Omega_B$.

As we have defined methods to handle all kinds of constraints, we will now explain the main approach of this work.
The proposed GSDO algorithm is implemented in the following three stages: (1) find initial feasible point, (2) generate multiple feasible points and (3) find global optimal point.

\subsection{Stage-1: Find Initial Feasible Point}
To start the procedure for obtaining the global solution to the problem (\ref{prob:originalProb}), an initial feasible point is required.
Accordingly, we solve an optimization problem modeled by surrogate constraints instead of the original constraints.
Let $m'$ be the number of quantifiable constraints in the problem (\ref{prob:originalProb}).
Here, each of the original constraint function $g_j$ from the problem (\ref{prob:originalProb}), if quantifiable, is substituted by its corresponding surrogate function $g_j^s:\mathbb{R}^d \rightarrow \mathbb{R}$ using Equation (\ref{rbf_cubic}).
Note that surrogate function cannot be created for nonquantifiable constraints (NRSK, NUSK or NUSH) due to lack of sufficient information from the constraint, and thus, need to be addressed in a different manner.
We now formulate the following new optimization problem which handles both quantifiable and nonquantifiable constraints:
\begin{equation}
	\tag{P2}
	\label{prob:feasiblePt}
	\begin{aligned}
		\max \, \, &z\\
		\text{ st. } &g_j^s(x) - z \geq 0 \quad \text{ for } j=1,\ldots,m'\\
		&x \geq l + ze\\
		&x \leq u - ze\\
		&z \geq 0\\
		& g_{c}(x) \geq 0
	\end{aligned}
\end{equation}
where $x \in \mathbb{R}^d$, $z \in \mathbb{R}$ and $e \in \mathbb{R}^d$ such that $e_i = 1$ for $i=1,\ldots,d$.

The first set of constraints using surrogate approximates corresponds to the actual nonlinear constraints.
Here, only quantifiable constraints (QRSK and QUSK) from problem (\ref{prob:originalProb}) are taken into account.
The nonquantifiable constraints are indirectly handled via the last constraint $g_c$.
In case problem (\ref{prob:feasiblePt}) is feasible, we denote its optimal solution as $(x^*,z^*)$.
Maximization of variable $z$ ensures that the optimal solution $x^*$ is strictly feasible for problem (\ref{prob:feasiblePt}) if $g_c(x^*) > 0$.
If surrogate constraint $g_j^s$ is a very close approximation to $g_j$ for all $j=1,\ldots,m'$, and $x^*$ is feasible for the nonquantifiable constraints, then heuristically, $x^*$ has a high chance of being feasible for the problem (\ref{prob:originalProb}).
\begin{algorithm}[h]
	\caption{Find Feasible Point}
	\label{algo:feasible}
	\begin{algorithmic}[1]
		\State Initialize: Set number of evaluations $k=0$.
		\While {Algorithm \ref{algo:rbf_rank} returns FALSE and $k \leq T_{max}$}
			\State Create a set of $T_{LH}$ points i.e. $Y=\{y_1,\ldots, y_{T_{LH}}\}$ inside $\Omega_B$ using LH search such that $|T_{LH}| <= T_{max}-k$.
			\State Evaluate $f$ and $\mathcal{G}$ for each point in $Y$.
			\State Filter each point in $Y$ using Algorithm (\ref{algo:filter}) and update $k \leftarrow k+T_{LH}$.
			\State Run Algorithm \ref{algo:rbf_rank} to check if sufficient points are available for RBF fit.
		\EndWhile
		\If {a feasible point is found}
			\State \textbf{Terminate}.
		\Else
			\State Set $k \leftarrow \text{col}(\mathcal{X})$.
			\While{$k \leq T_{max}$}
				\State Fit surrogate models for QRSK and QUSK constraints using Algorithm (\ref{algo:surrogate_fit}).
				\State Solve the problem (\ref{prob:feasiblePt}). Compute constraint $g_c$ using Algorithm (\ref{algo:knn}).
				\If {solution $x^*$ exists}
				\State Perform constraint evaluations at $x^*$ and update $k \leftarrow k+1$.
				\State Filter $x^*$ using Algorithm (\ref{algo:filter}).
				\If {$x^* \in \Omega$} 
					\State \textbf{Terminate}.
				\EndIf
				\Else
					\State Create a set $T$ of $t=\text{min }\{d+1, T_{max}-k\}$ new points using LH search.
					\For {each point in $T$}
						\State Evaluate $f$ and $\mathcal{G}$ and filter it using Algorithm (\ref{algo:filter}).
						\State If $\text{col}(\mathcal{X}_F) > 0$ then \textbf{terminate}.
					\EndFor
					\State Update $k \leftarrow k+|T|$.
				\EndIf
			\EndWhile
		\EndIf
	\end{algorithmic}
\end{algorithm}

The working principle of Stage-1 is summarized in Algorithm (\ref{algo:feasible}).
The initial set of points required for modeling the surrogate functions are obtained using Latin Hypercube algorithm (\cite{forrester2008engineering}).
The objective function and constraints from problem (\ref{prob:originalProb}) are computed over these points.
Next, Algorithm (\ref{algo:filter}) is applied to segregate these evaluated points into mutually exclusive sets $\mathcal{X}_F$, $\mathcal{X}_I$, $\mathcal{X}_S$, $\mathcal{X}_U$ and $\mathcal{X}_H$.
Subsequently, for each quantifiable constraint, Algorithm (\ref{algo:rbf_rank}) is applied to check whether the required conditions for RBF surrogate model (Equation (\ref{eq:rbf_linear_system})) are met.
If true, RBF surrogates are constructed using Algorithm (\ref{algo:surrogate_fit}); otherwise, more points are explored using Latin Hypercube and procedure is repeated till some termination criterion is met.
In the case of nonquantifiable constraints, the computation of $g_c$ constraint of problem (\ref{prob:feasiblePt}) is performed using Algorithm (\ref{algo:knn}), which returns a positive or negative value based on the classification of a point.
Consequently, if a solution $x \in \mathbb{R}^n$ for the problem (\ref{prob:feasiblePt}) exists, then it is likely to be distant and separated from the cluster of points in the sets $\mathcal{X}_S$, $\mathcal{X}_U$ and $\mathcal{X}_H$ and closer to the points in the sets $\mathcal{X}_F$ or $\mathcal{X}_I$.
Effectively, this leads to higher possibility of $x$ being more feasible to the original DFO problem when nonquantifiable constraints are present.
Next, $x$ is evaluated for the problem (\ref{prob:originalProb}).
If $x$ is found feasible, the Algorithm (\ref{algo:feasible}) is terminated; otherwise the constraints $\mathcal{G}^s$ and $g_c$ are remodeled using it.
However, if no solution exists for the problem (\ref{prob:feasiblePt}), then new points are generated through Latin Hypercube search and evaluated for the original constraints.
Note that if any of these points is found feasible, then Algorithm (\ref{algo:feasible}) is terminated.
Otherwise, these points are filtered into their respective sets using Algorithm (\ref{algo:filter}) and the constraints $\mathcal{G}^s$ and $g_c$ are remodeled.
The problem (\ref{prob:feasiblePt}) with remodeled constraints is solved again and the above process is repeated either till a feasible point is found for the original problem or some termination criterion is met.

\subsection {Stage-2: Generate multiple feasible points}
\begin{algorithm}[h]
	\caption {Generate Multiple Feasible Points}
	\label{algo:multiple}
	\begin{algorithmic}[1]
		\State Initialization: Set $\eta_f \leftarrow \text{col}(\mathcal{X}_F)$. 
		\While {$\eta_f < \eta_{max}$ and $\text{col}(\mathcal{X}) \leq T_{max}$}
		\State Solve the problem (\ref{eq:multiple_feasible_opt}). Constraint $g_c$ is computed using Algorithm (\ref{algo:knn}).
		\State Let $x$ be the solution (if exists) of problem (\ref{eq:multiple_feasible_opt}).
		\If {Solution does not exist}
			\State Generate $x$ using Algorithm (\ref{algo:random_pt}).
		\EndIf
		\State Evaluate $f$ and $\mathcal{G}$ at $x$.
		\State Filter $x$ using Algorithm (\ref{algo:filter}).
		\If {$x$ is feasible}
			\State Set $\eta_f \leftarrow \eta_f+1$. 
		\EndIf
		\State Fit surrogate models for QRSK and QUSK constraints using Algorithm (\ref{algo:surrogate_fit}).
	\EndWhile
	\end{algorithmic}
\end{algorithm}
The next part of the proposed approach is to find multiple points which are well spread within the feasible region.
This step facilitates generation of better surrogate models due to increased exposure to the feasible region.
Accordingly, we formulate an optimization problem in which a new point is computed by maximizing its distance from the feasible points in the set $\mathcal{X}_F$ while maintaining the feasibility constraints.
The optimization problem is given as:
\begin{equation}
	\tag{P3}
	\label{eq:multiple_feasible_opt}
	\begin{aligned}
		\max \, \, & y\\
		\text{ st. } &g_j^s(x)  \geq 0 \quad \text{ for } j=1,\ldots,m'\\
		&||x - x^i ||_2^2 \geq y^2 \text{ for } i=1,\ldots,|\mathcal{X}_F|\\
		&x \in \Omega_B\\
		& g_{c}(x) \geq 0\\
	\end{aligned}
\end{equation}
where $x \in \mathbb{R}^d$, $y \in \mathbb{R}_+$ and $x^i \in \mathcal{X}_F$.
Note that the problem (\ref{eq:multiple_feasible_opt}) is always feasible, as any $x^i \in \mathcal{X}_F$ will always satisfy the given constraints.
If the problem (\ref{eq:multiple_feasible_opt}) has an optimal solution $(x^*,y^*) \in \mathbb{R}^{d+1}$ such that $x^* \notin \mathcal{X}_F$, then we have $y^* > 0$.
In the subsequent analysis, we denote $y^*$ with $\Delta$, which indicates the distance of the solution $x^*$ from the points in the set $\mathcal{X}_F$.
Due to positive value of $\Delta$, the above formulation leads to exploration of the feasible space.
The working principle of Stage-2 is summarized in Algorithm (\ref{algo:multiple}).
If the optimization problem (\ref{eq:multiple_feasible_opt}) fails to obtain a solution $x^*$ such that $x^* \notin \mathcal{X}_F$, we create a random point in the vicinity of last feasible point from set $\mathcal{X}_F$ using Algorithm (\ref{algo:random_pt}).
This ensures evaluation of a new point in each iteration of Algorithm (\ref{algo:multiple}).

\subsection{Stage-3: Find Global Optimal Point}
\begin{algorithm}[htp]
	\caption{Find Global Optimal Point}
	\label{algo:global}
	\begin{algorithmic}[1]
	\State  Initialize: Set number of evaluations $k=0$.
	\While {$k < T_{max}$}
		\State Fit the surrogate model for objective function using points from $\mathcal{X}_F\cup \mathcal{X}_I\cup \mathcal{X}_S$.
		\State Fit the surrogate models for QRSK and QUSK constraints using Algorithm (\ref{algo:surrogate_fit}).
	\State Solve the problem (\ref{eq:local_opt}) and find multiple local optima.
	\If {solution exists}
		\State Create the set $\mathcal{X}_{opt}$ by sorting the solutions in ascending order based on their function values.
	\EndIf
	\If {solution does not exist or $\mathcal{X}_{opt} \subset \mathcal{X}$}
		\State Set the search type to ``Exploration".
	\Else
		\State Set the search type based on some random number $t$.
	\EndIf
	\If {$k > K_{Global}$ and $t < C_g$ }
		\State Set the search type to ``Exploitation" and choose the first point $x \in \mathcal{X}_{opt}$ such that $x \notin \mathcal{X}$.
	\Else
		\State Set the search type to ``Exploration". Solve the problem (\ref{eq:multiple_feasible_opt}) to obtain $\Delta$.
		\If {$\Delta = 0$}
			\State Set $\Delta = \min \{\delta_{lu}, 1\}$.
		\EndIf
		\If {$\Delta < \Delta_{min}$}
			\State \textbf{Terminate} the algorithm.
		\Else
			\State Solve the problem (\ref{eq:global_opt}) to find the global optimal solution $x$. 
			\If {problem (\ref{eq:global_opt}) is infeasible}
				\State Create a random point $x$ using Algorithm (\ref{algo:random_pt}).
			\EndIf
		\EndIf
	\EndIf
	\State Evaluate $f$ and $\mathcal{G}$ at $x$. Update $k \leftarrow k+1$.
	\State Filter $x$ using Algorithm (\ref{algo:filter}).
	\EndWhile
	\end{algorithmic}
\end{algorithm}

The third and final stage of GSDO describes the methodology to solve problem (\ref{prob:originalProb}) in a global sense. 
This stage involves a combination of exploitation and exploration steps.

The exploitation step is performed by solving the following optimization problem: 
\begin{equation}
	\tag{P4}
	\label{eq:local_opt}
	\begin{aligned}
		\min \, \, & s(x)\\
		\text{ st. } &g_j^s(x)  \geq 0 \quad \text{ for } j=1,\ldots,m'\\
		&x \in \Omega_B\\
		& g_{c}(x) \geq 0
	\end{aligned}
\end{equation}
where $s:\mathbb{R}^d\rightarrow \mathbb{R}$ is a surrogate approximation to the objective function of problem (\ref{prob:originalProb}).
This problem is solved by a multistart strategy using different starting points, which leads to the creation of solution pool containing different local minima of the surrogate model.
The points in the solution pool are sorted in ascending order of their surrogate objective function values and stored in the set $\mathcal{X}_{opt}$.
Note that the above strategy can sometimes produce at most one solution i.e. $|\mathcal{X}_{opt}|=1$.
In exploitation step, the first point $x$ from $\mathcal{X}_{opt}$ such that $x \notin \mathcal{X}$ is chosen for the evaluation of objective function and constraints. 
This ensures exploitation at the best solution obtained from the surrogate model.
In case $\mathcal{X}_{opt}=\emptyset$ or $\mathcal{X}_{opt} \subset \mathcal{X}$, the solver invokes the exploration step.

The exploration search is carried out by solving the following optimization problem:
\begin{equation}
	\tag{P5}
	\label{eq:global_opt} 
	\begin{aligned}
		\min \, \, & s(x)\\
		\text{ st. } &g_j^s(x)  \geq 0 \quad \text{ for } j=1,\ldots,m'\\
		&||x - x^i ||_2^2 \geq \Delta^2 \text{ for } i=1,\ldots,|\mathcal{X}_F|\\
		&x \in \Omega_B\\
		& g_{c}(x) \geq 0
	\end{aligned}
\end{equation}
where $x\in \mathbb{R}^d$ and $x^i \in \mathcal{X}_F$.
The computation of constraint $g_c$ in the problems \ref{eq:multiple_feasible_opt}, \ref{eq:local_opt} and \ref{eq:global_opt} is done through Algorithm (\ref{algo:knn}).
Note that $\Delta$ is obtained by solving the problem (\ref{eq:multiple_feasible_opt}) and is always non-negative.
Further, if $\Delta = 0$, it is set to $\min\{\delta_{lu},1\}$ in Algorithm (\ref{algo:global}) in order to always make it positive.
Consequently, the second constraint in the problem (\ref{eq:global_opt}) ensures better exploration of the feasible space by enforcing a positive distance (greater than or equal to $\Delta$) between the corresponding optimal solution and the points in the set $\mathcal{X}_F$.

Algorithm (\ref{algo:global}) mainly works by switching between the exploration and exploitation steps using a random number, which is generated in each iteration.
Once the type of step is determined, either problem (\ref{eq:local_opt}) or problem (\ref{eq:global_opt}) is solved accordingly.
During any iteration, if the exploitation step fails because of a repeated solution or infeasibility of the problem (\ref{eq:local_opt}), we switch to the exploration step.
However, in case the exploration step fails, a random point is generated using Algorithm (\ref{algo:random_pt}).
Subsequently, the expensive objective function and constraints are evaluated at the solution, and the surrogate models are updated for the next iteration.
The points generated by the Algorithm (\ref{algo:global}) in each iteration enrich the required sets for K-Nearest Neighbor method, thus leading to better classification with $g_c$ constraint.
Also, if these points are filtered into the sets $\mathcal{X}_F$, $\mathcal{X}_I$, $\mathcal{X}_S$ or $\mathcal{X}_U$, they are helpful in creating better surrogate models for QRSK and QUSK constraints.
Thus, heuristically, the likelihood of finding feasible points for problem (\ref{eq:global_opt}) increases with the number of iterations.
Overall, as the iterations proceed, the feasible search space $\Omega$ is gradually populated, effectively leading to higher possibility of finding the global optimal solution.

However, for practical purposes, a criterion is needed for termination of the algorithm within finite iterations. 
Accordingly, a lower bound $\Delta_{min}$ on the distance $\Delta$ is chosen as the termination criterion in addition to maximum allowed budget.
The working principle of Stage-3 is summarized in Algorithm \ref{algo:global}.

\section{Numerical Experiments} \label{numerical}
The entire algorithm was implemented in C language.
We tested our proposed approach on a diverse collection of test problems taken from \cite{regis2014constrained}.
We considered 18 test problems, out of which 4 are engineering application problems.
Table \ref{table:probleminfo} shows the list of test problems considered, along with information about their corresponding dimension, number of constraints involved, and the global optimum (or the best known solution).
The test problems possess a wide diversity in terms of dimension and number of constraints.
The dimensions vary from 2 to 20 while the number of constraints vary from 1 to 13, which is a descent set for testing constrained DFO problems.
We further divide the test problems into four different sets, where the first set has all the 18 problems while the remaining three sets have 16 problems each.
These later sets consider only those problems where the number of constraints are greater than 1.
So G3MOD and GTCD test problems, which involve only one constraint, are excluded from these sets.
We now define these four sets as:
\begin{itemize}
	\item Set-1: The constraints in all the problems are treated as QRSK.
	\item Set-2: The first constraint in each problem is treated as NRSK and the remaining constraints as QRSK.
	\item Set-3: The first constraint in each problem is treated as QUSK and the remaining constraints as QRSK.
	\item Set-4: The first constraint in each problem is treated as NUSK and the remaining constraints as QRSK.
\end{itemize}
\begin{table*}[t]
	\center
	\caption {Test Problems For Global Optimization}
	\begin{tabular}{|p{5cm}|c|c|c|}
		\hline
		Problem Name&Dimension&Constraints&Global Optimum\\
		\hline
		G1&13&9&\num{-15}\\
		G2&10&2&\num{-0.4}\\
		G3MOD&20&1&\num{-0.69}\\
		G4&5&6&\num{-30665.539}\\
		G5MOD&4&5&\num{5126.50}\\
		G6&2&2&\num{-6961.8139}\\
		G7&10&8&\num{24.3062}\\
		G8&2&2&\num{-0.0958}\\
		G9&7&4&\num{680.6301}\\
		G10&8&6&\num{7049.3307}\\
		G18&9&13&\num{-0.8660}\\
		G19&15&5&\num{32.6556}\\
		G24&2&2&\num{-5.5080}\\
		Gas transmission compressor design (GTCD)&4&1&\num{2964893.85}\\
		Hesse&6&6&\num{-310}\\
		Pressure vessel design (PVD)&4&3&\num{5804.45}\\
		Speed reducer (SR)&7&11&\num{2994.42}\\
		Welded beam (WB)&4&6&\num{1.7250}\\
		\hline
	\end{tabular}
	\label {table:probleminfo}
\end{table*}

\begin{table*}[t!]
	\center
	\caption{Median Optimal Values attained by Solvers for Set-1}
	\begin{tabular}{|c|c|c|c|c|c|c|c|c|}
		\hline
		Name &\multicolumn{2}{c|}{GSDO} &\multicolumn{2}{c|}{NOMAD} &\multicolumn{2}{c|}{SHGO} &\multicolumn{2}{c|}{DE}\\
		\hline
		&$N_s$&Value&$N_s$&Value&$N_s$&Value&$N_s$&Value\\
		\hline
		G1&30&\num{-1.328e+01}&17&\num{-6.000e+00}&30&\num{-6.000e+00}&0&-\\
		G2&30&\num{-7.755e-03}&30&\num{-1.018e-02}&0&-&30&\num{-8.627e-03}\\
		G3MOD&30&\num{-1.953e-01}&30&\num{-1.098e-04}&30&\num{-0.000e+00}&0&-\\
		G4&30&\num{-3.067e+04}&30&\num{-3.022e+04}&30&\num{-2.859e+04}&30&\num{-2.926e+04}\\
		G5MOD&30&\num{5.127e+03}&7&\num{6.295e+03}&30&\num{8.880e+03}&0&-\\
		G6&30&\num{-6.665e+03}&14&\num{-3.745e+03}&0&-&2&\num{-1.246e+03}\\
		G7&30&\num{4.319e+01}&20&\num{4.755e+02}&0&-&0&-\\
		G8&30&\num{-9.582e-02}&29&\num{-1.166e-02}&30&\num{1.024e-01}&9&\num{-1.901e-04}\\
		G9&30&\num{1.234e+03}&30&\num{1.582e+03}&0&-&18&\num{9.615e+03}\\
		G10&30&\num{8.482e+03}&3&\num{1.332e+04}&0&-&0&-\\
		G18&21&\num{-3.084e-01}&0&-&0&-&0&-\\
		G19&30&\num{2.181e+02}&28&\num{1.847e+03}&30&\num{9.125e+02}&0&-\\
		G24&30&\num{-5.504e+00}&30&\num{-5.235e+00}&30&\num{-4.250e+00}&30&\num{-4.753e+00}\\
		GTCD&30&\num{3.337e+06}&30&\num{5.636e+06}&30&\num{1.206e+07}&30&\num{6.963e+06}\\
		Hesse&30&\num{-3.100e+02}&30&\num{-1.841e+02}&30&\num{-1.320e+02}&30&\num{-1.583e+02}\\
		PVD&30&\num{5.900e+03}&30&\num{9.481e+03}&30&\num{7.329e+03}&14&\num{8.012e+03}\\
		SR&30&\num{2.994e+03}&30&\num{2.983e+03}&30&\num{3.089e+03}&4&\num{4.208e+03}\\
		WB&29&\num{3.063e+00}&11&\num{7.244e+00}&0&-&0&-\\
		\hline
	\end{tabular}
	\label{table:median_qrsk}
\end{table*}

\begin{table*}[t!]
	\center
	\caption{Median Optimal Values attained by Solvers for Set-2}
	\begin{tabular}{|c|c|c|c|c|c|c|c|c|}
		\hline
		Name &\multicolumn{2}{c|}{GSDO} &\multicolumn{2}{c|}{NOMAD} &\multicolumn{2}{c|}{SHGO} &\multicolumn{2}{c|}{DE}\\
		\hline
		&$N_s$&Value&$N_s$&Value&$N_s$&Value&$N_s$&Value\\
		\hline
		G1&30&\num{-1.330e+01}&27&\num{-1.260e+01}&30&\num{-6.000e+00}&0&-\\
		G2&30&\num{-3.961e-02}&30&\num{-1.747e-02}&0&-&30&\num{-1.281e-02}\\
		G4&30&\num{-3.067e+04}&30&\num{-3.049e+04}&30&\num{-2.859e+04}&30&\num{-2.964e+04}\\
		G5MOD&29&\num{5.127e+03}&17&\num{5.607e+03}&30&\num{8.880e+03}&9&\num{7.231e+03}\\
		G6&0&-&30&\num{-4.206e+03}&0&-&1&\num{-3.207e+03}\\
		G7&30&\num{4.673e+01}&30&\num{8.464e+01}&0&-&0&-\\
		G8&29&\num{-8.640e-02}&29&\num{-8.772e-02}&30&\num{1.024e-01}&24&\num{-5.740e-03}\\
		G9&30&\num{9.159e+02}&24&\num{8.830e+02}&30&\num{1.183e+03}&29&\num{9.688e+03}\\
		G10&20&\num{9.857e+03}&0&-&0&-&0&-\\
		G18&0&-&7&\num{-3.755e-01}&0&-&0&-\\
		G19&30&\num{1.443e+02}&29&\num{1.048e+03}&30&\num{9.125e+02}&1&\num{1.661e+04}\\
		G24&30&\num{-5.297e+00}&30&\num{-5.450e+00}&30&\num{-4.750e+00}&30&\num{-5.045e+00}\\
		Hesse&30&\num{-3.100e+02}&30&\num{-2.651e+02}&30&\num{-1.320e+02}&30&\num{-1.937e+02}\\
		PVD&28&\num{6.717e+03}&30&\num{7.518e+03}&30&\num{7.329e+03}&21&\num{7.802e+03}\\
		SR&30&\num{2.994e+03}&30&\num{2.933e+03}&30&\num{3.089e+03}&21&\num{3.778e+03}\\
		WB&29&\num{4.883e+00}&21&\num{5.866e+00}&0&-&1&\num{4.317e+00}\\
		\hline
	\end{tabular}
	\label{table:median_nrsk}
\end{table*}

\begin{table*}[t!]
	\center
	\caption{Median Optimal Values attained by Solvers for Set-3}
	\begin{tabular}{|c|c|c|c|c|c|c|c|c|}
		\hline
		Name &\multicolumn{2}{c|}{GSDO} &\multicolumn{2}{c|}{NOMAD} &\multicolumn{2}{c|}{SHGO} &\multicolumn{2}{c|}{DE}\\
		\hline
		&$N_s$&Value&$N_s$&Value&$N_s$&Value&$N_s$&Value\\
		\hline
		G1&0&-&30&\num{-1.172e+01}&30&\num{-6.000e+00}&0&-\\
		G2&30&\num{-2.974e-02}&30&\num{-1.894e-02}&0&-&30&\num{-1.141e-02}\\
		G4&30&\num{-3.067e+04}&30&\num{-3.048e+04}&30&\num{-2.859e+04}&30&\num{-2.951e+04}\\
		G5MOD&30&\num{5.127e+03}&28&\num{5.817e+03}&30&\num{8.880e+03}&9&\num{6.850e+03}\\
		G6&30&\num{-6.959e+03}&26&\num{-4.206e+03}&0&-&1&\num{-2.745e+03}\\
		G7&30&\num{4.186e+01}&30&\num{1.103e+02}&0&-&0&-\\
		G8&30&\num{-9.582e-02}&30&\num{-9.211e-02}&30&\num{1.024e-01}&27&\num{-6.989e-03}\\
		G9&30&\num{9.276e+02}&29&\num{8.772e+02}&30&\num{1.183e+03}&28&\num{1.028e+04}\\
		G10&29&\num{9.183e+03}&26&\num{9.000e+03}&0&-&0&-\\
		G18&0&-&1&\num{-8.500e-02}&0&-&0&-\\
		G19&30&\num{1.416e+02}&30&\num{4.668e+02}&30&\num{9.125e+02}&0&-\\
		G24&30&\num{-5.508e+00}&30&\num{-5.400e+00}&30&\num{-4.750e+00}&30&\num{-5.064e+00}\\
		Hesse&30&\num{-3.100e+02}&30&\num{-2.790e+02}&30&\num{-1.320e+02}&30&\num{-1.846e+02}\\
		PVD&29&\num{7.762e+03}&30&\num{7.100e+03}&30&\num{7.329e+03}&25&\num{7.885e+03}\\
		SR&30&\num{2.994e+03}&30&\num{2.938e+03}&30&\num{3.089e+03}&22&\num{3.820e+03}\\
		WB&30&\num{2.920e+00}&20&\num{5.028e+00}&0&-&2&\num{6.811e+00}\\
		\hline
	\end{tabular}
	\label{table:median_qusk}
\end{table*}

\begin{table*}[t!]
	\center
	\caption{Median Optimal Values attained by Solvers for Set-4}
	\begin{tabular}{|c|c|c|c|c|c|c|c|c|}
		\hline
		Name &\multicolumn{2}{c|}{GSDO} &\multicolumn{2}{c|}{NOMAD} &\multicolumn{2}{c|}{SHGO} &\multicolumn{2}{c|}{DE}\\
		\hline
		&$N_s$&Value&$N_s$&Value&$N_s$&Value&$N_s$&Value\\
		\hline
		G1&0&-&2&\num{-1.160e+01}&30&\num{-6.000e+00}&0&-\\
		G2&30&\num{-2.109e-02}&30&\num{-1.946e-02}&0&-&30&\num{-1.125e-02}\\
		G4&30&\num{-3.067e+04}&30&\num{-3.055e+04}&30&\num{-2.859e+04}&30&\num{-2.965e+04}\\
		G5MOD&30&\num{5.127e+03}&27&\num{5.976e+03}&30&\num{8.880e+03}&11&\num{6.840e+03}\\
		G6&0&-&28&\num{-4.206e+03}&0&-&0&-\\
		G7&30&\num{4.600e+01}&29&\num{9.296e+01}&0&-&0&-\\
		G8&30&\num{-8.108e-02}&18&\num{-6.936e-02}&30&\num{1.024e-01}&23&\num{-1.422e-03}\\
		G9&30&\num{9.734e+02}&18&\num{8.607e+02}&30&\num{1.183e+03}&28&\num{4.089e+03}\\
		G10&25&\num{1.483e+04}&3&\num{1.666e+04}&0&-&0&-\\
		G18&0&-&2&\num{-7.882e-02}&0&-&0&-\\
		G19&30&\num{1.512e+02}&30&\num{9.532e+02}&30&\num{9.125e+02}&1&\num{1.154e+04}\\
		G24&30&\num{-5.334e+00}&30&\num{-5.435e+00}&30&\num{-4.750e+00}&30&\num{-4.836e+00}\\
		Hesse&30&\num{-3.100e+02}&30&\num{-2.842e+02}&30&\num{-1.320e+02}&30&\num{-1.868e+02}\\
		PVD&30&\num{7.407e+03}&30&\num{7.315e+03}&30&\num{7.329e+03}&23&\num{7.723e+03}\\
		SR&30&\num{2.994e+03}&30&\num{2.937e+03}&30&\num{3.089e+03}&16&\num{4.032e+03}\\
		WB&26&\num{5.047e+00}&20&\num{4.923e+00}&0&-&5&\num{9.138e+00}\\
		\hline
	\end{tabular}
	\label{table:median_nusk}
\end{table*}

\subsection{Algorithms for Comparison}
There are not many open source solvers which are designed specifically to address constrained DFO in a global sense.
Hence, we compared our approach with generic but state of the art solvers like Nonlinear Optimization by Mesh Adaptive Direct Search (NOMAD), Simplicial Homology Global Optimization (SHGO) and differential evolution (DE). 
NOMAD by \cite{Le2011a} is a mesh adaptive direct search algorithm to solve most variants of blackbox optimization problems. 
SHGO algorithm by \cite{endres2018simplicial} is a global optimization approach that utilizes the homological properties of the objective function surface.
DE algorithm by \cite{storn1997differential} is a stochastic, population based metaheuristic approach designed to solve global optimization problems without the use of gradients.
This approach, however, requires large number of function evaluations.
The SCIPY implementation of SHGO and DE were considered for the comparison tests.

\subsection{Setup of different algorithms for numerical experiments}
A simulation of the test problem consists of computing the expensive objective function and evaluation of all its constraints.
For all problems in Set-1, the maximum budget for simulations was set to $15(d+1)$, where $d$ is the dimension of the test problem.
For the remaining sets, which are more difficult to solve, the maximum budget for simulations was increased to $30(d+1)$.
The reason for choosing such a small number is that finding a global optimum for these highly expensive DFO problems under a low budget restriction is an important practical concern.
For NOMAD, the variable neighborhood search (VNS) parameter was enabled for solving the problem in a global sense, while the rest of the parameters were set to their respective defaults.
The SHGO and DE solvers were executed such that the respective algorithms were terminated as soon as the maximum limit of function evaluations was reached.
The parameters for both these solvers were set to their respective defaults.
For GSDO, we set the parameters as follows:
$T_{max} = 15(d+1)$ for Set-1 and $30(d+1)$ for remaining sets, $c_1=1$, $c_2=-1$, $c_3=-10$, $c_4=-100$, $\delta_r=10$, $\delta_d=100$, $\eta_{max}= d+1$, $C_g = 0.50$, $K_{Global} = d+1$ and $\Delta_{min} = 10^{-5}$.
We used differential evolution algorithm (\cite{lampinen2002constraint}) for solving the surrogate optimization problems (\ref{prob:feasiblePt}) in Stage-1, (\ref{eq:multiple_feasible_opt}) in Stage-2, and (\ref{eq:local_opt}) and (\ref{eq:global_opt}) in Stage-3.

\subsection{Results}
Tables \ref{table:median_qrsk} to \ref{table:median_nusk} show the optimal values attained by the four solvers under different constraint conditions.
The column header ``$N_s$" under each solver represents the number of times a particular solver was able to find a feasible solution out of 30 trials with different starting points.
We define a trial as successful for the solver, if at least one feasible solution was found out of 30 trials.
The median of optimal values attained by the solver over successful trials is reported under column header ``Value".
The blank entry indicated by ``$-$" shows that the corresponding solver was unsuccessful in finding any feasible solution in all the trials.

We observe that NOMAD, SHGO and DE were able to find at least one feasible solution in 17, 11 and 10 test cases in Set-1; 15, 10 and 12 test cases in Set-2; 16, 10 and 11 test cases in Set-3; and 16, 10 and 11 test cases in Set-4.
However, if we consider at least 15 successful trials out of 30 for each test case, then we see different trends.
NOMAD, SHGO and DE were successful in 13, 11 and 6 test cases in Set-1; 14,10 and 8 test cases in Set-2; 15, 10 and 8 test cases in Set-3; and 13, 10 and 8 test cases in Set-4.
Since, SHGO is agnostic to the starting point, all trials gave the same solution for any particular problem.
Our proposed approach, GSDO was successful on all the problems in Set-1, 14 test cases in Set-2 and Set-3, and 13 test cases in Set-4.
In addition, it can be noted that the proposed method was successful for at least 20 trials in each case.
Clearly, GSDO demonstrates good and consistent performance on all the problems.

\begin{figure}[t!]
	\subfloat[]{ \includegraphics[width=0.4\textwidth]{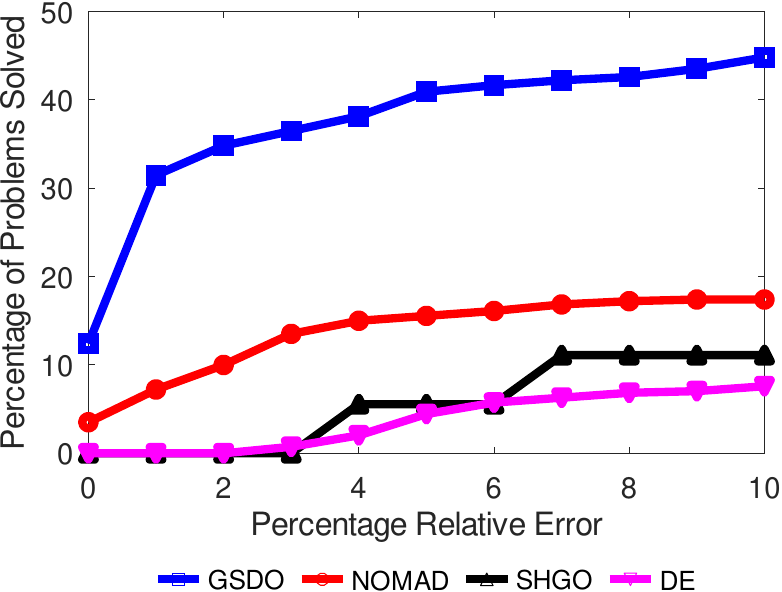}}
	\subfloat[]{ \includegraphics[width=0.4\textwidth]{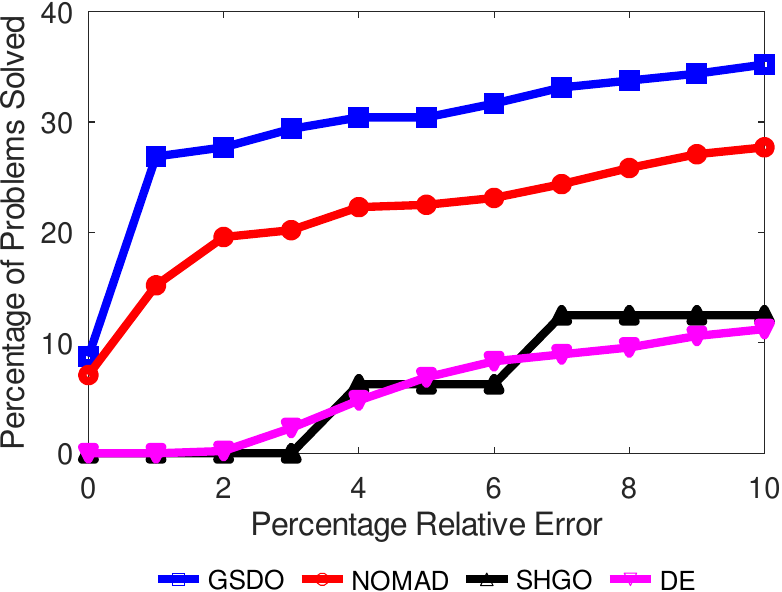}}
	\newline
	\subfloat[]{ \includegraphics[width=0.4\textwidth]{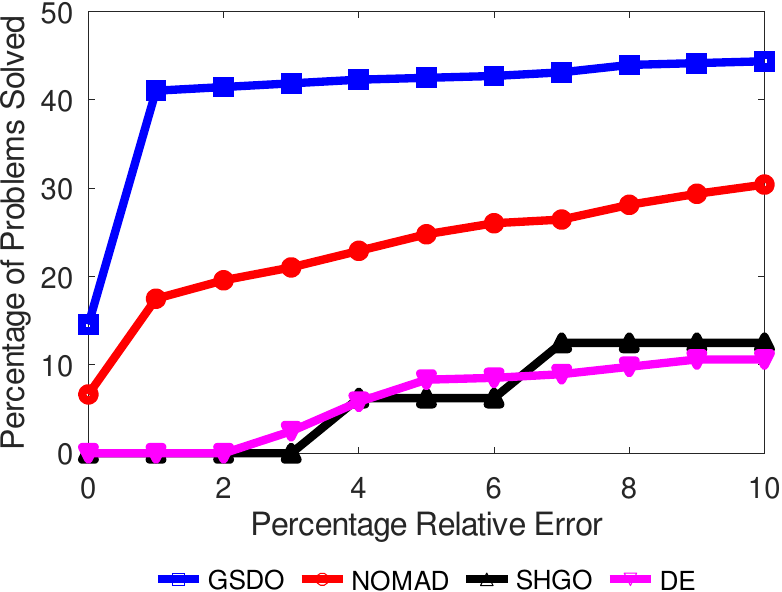}}
	\subfloat[]{ \includegraphics[width=0.4\textwidth]{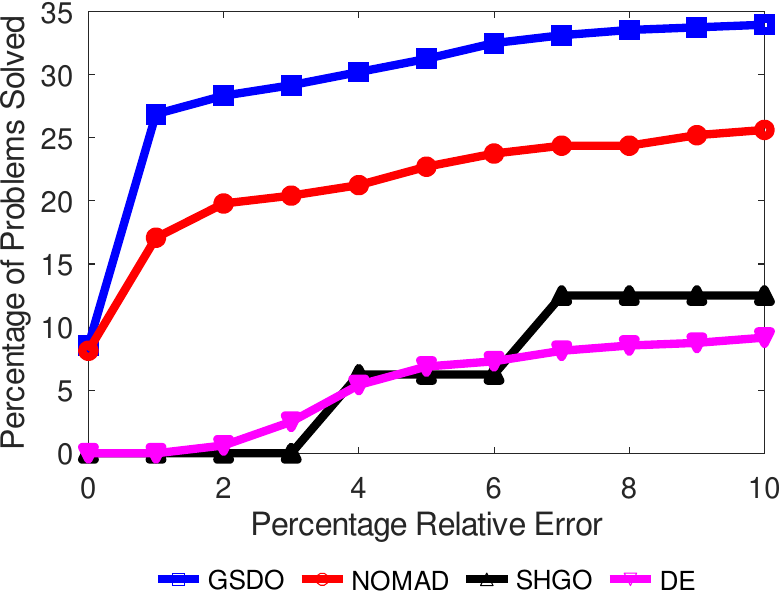}}
	\caption{Percentage of problems solved within percentage relative error (a) Set-1 (b) Set-2 (c) Set-3 (d) Set-4}
	\label{fig:errorplot}
\end{figure}

\begin{figure}[ht]
	\subfloat[]{ \includegraphics[width=0.4\textwidth]{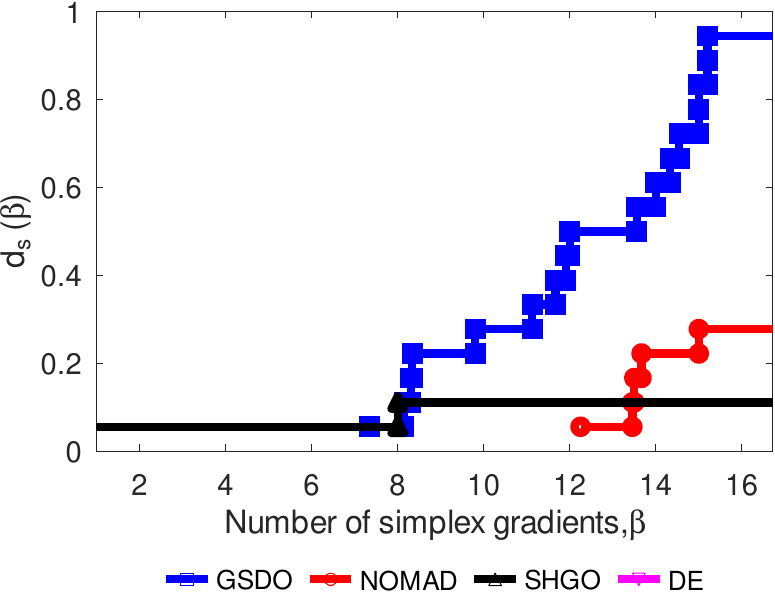}}
	\subfloat[]{ \includegraphics[width=0.4\textwidth]{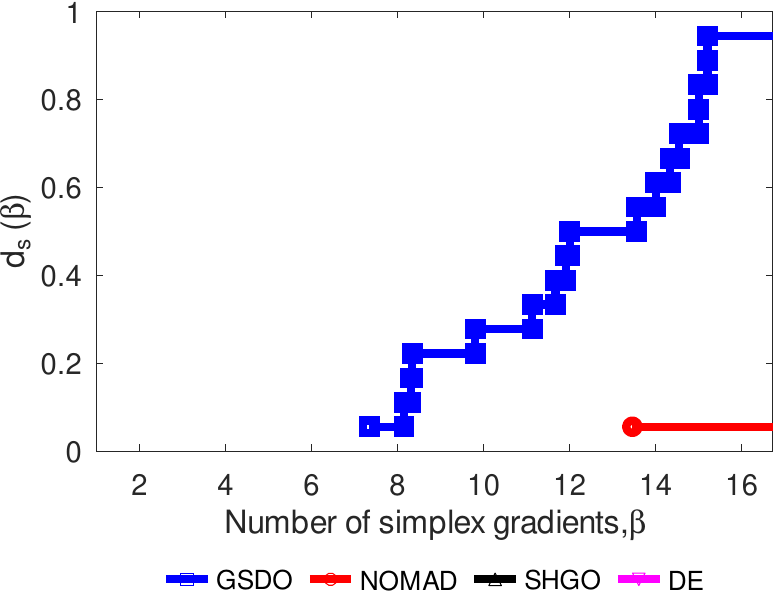}}
	\newline
	\subfloat[]{ \includegraphics[width=0.4\textwidth]{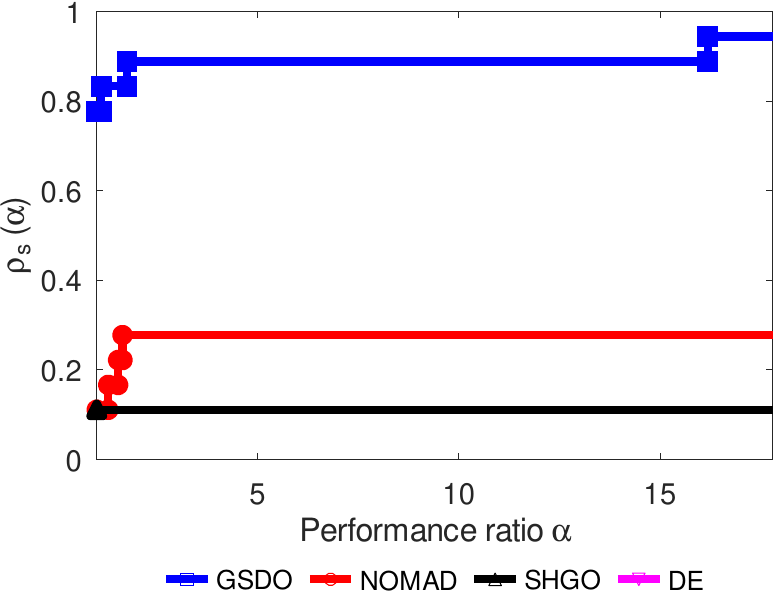}}
	\subfloat[]{ \includegraphics[width=0.4\textwidth]{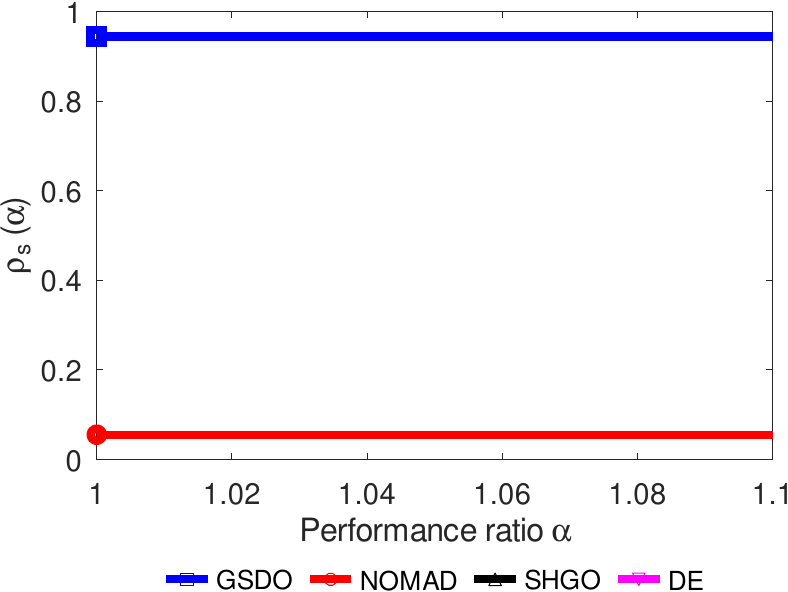}}
	\caption{Set-1: (a) Data profiles at $\tau=0.1$ (b) Data profiles at $\tau=0.01$ (c) Performance profiles at $\tau=0.1$ (d) Performance profiles at $\tau=0.01$}
	\label{fig:dataperfprofile_qrsk}
\end{figure}

\begin{figure}[h]
	\subfloat[]{ \includegraphics[width=0.4\textwidth]{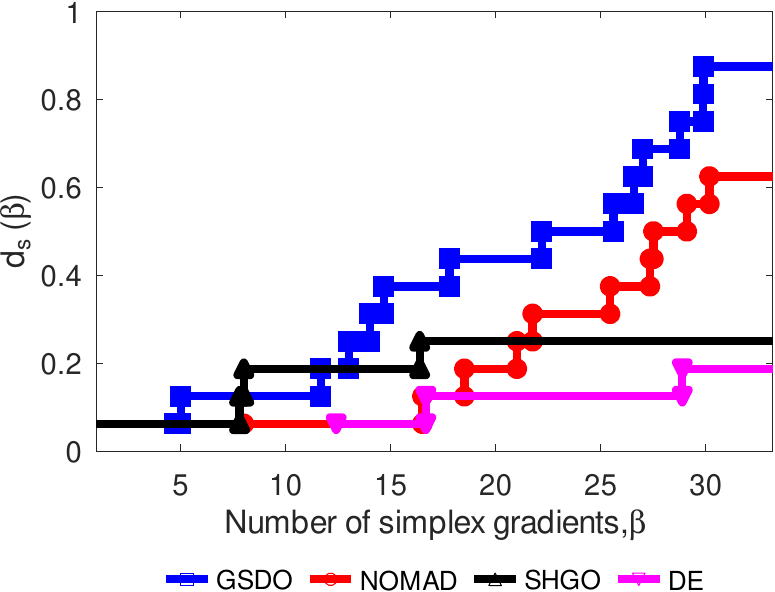}}
	\subfloat[]{ \includegraphics[width=0.4\textwidth]{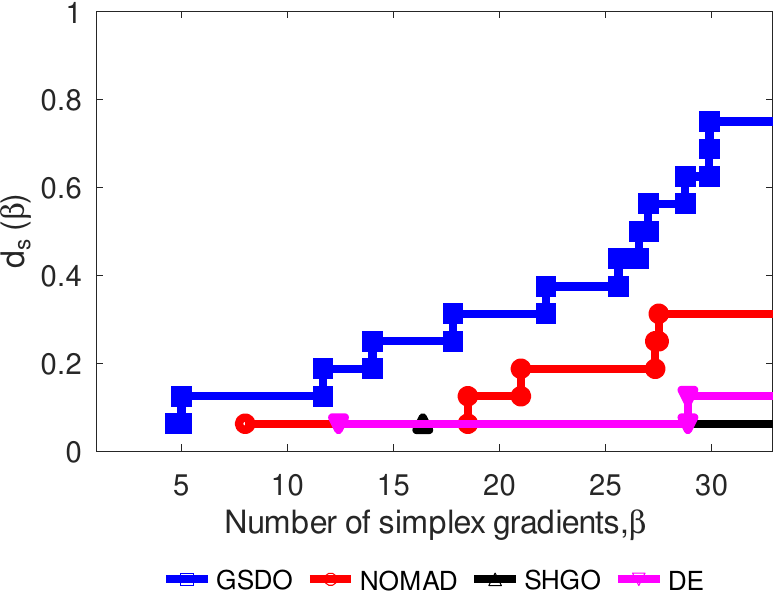}}
	\newline
	\subfloat[]{ \includegraphics[width=0.4\textwidth]{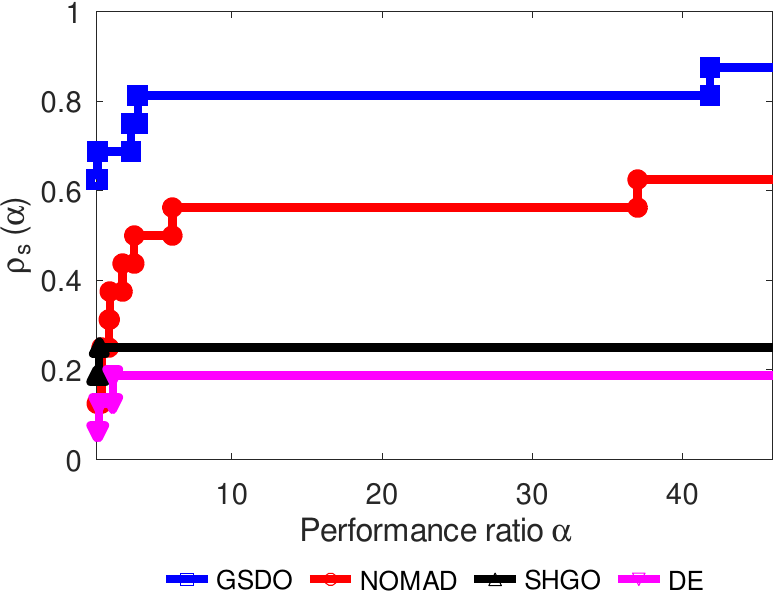}}
	\subfloat[]{ \includegraphics[width=0.4\textwidth]{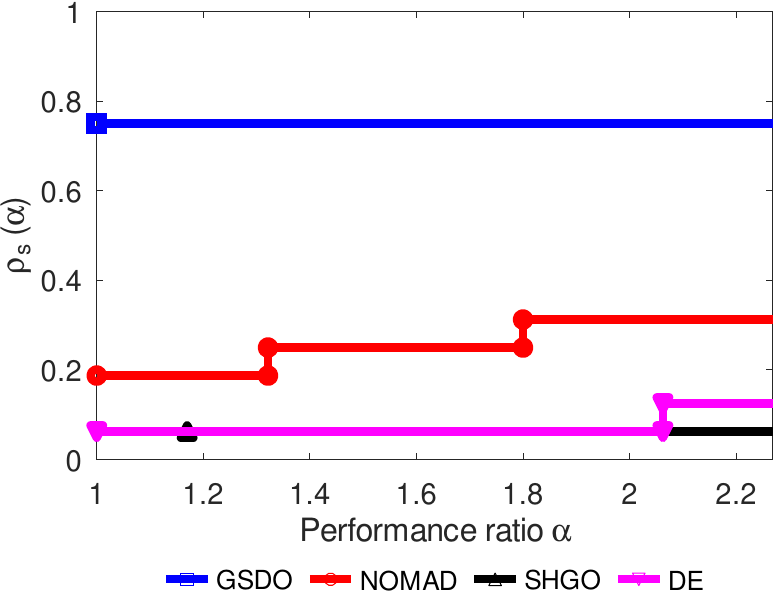}}
	\caption{Set-2: (a) Data profiles at $\tau=0.1$ (b) Data profiles at $\tau=0.01$ (c) Performance profiles at $\tau=0.1$ (d) Performance profiles at $\tau=0.01$}
	\label{fig:dataperfprofile_nrsk}
\end{figure}

\begin{figure}[ht]
	\subfloat[]{ \includegraphics[width=0.4\textwidth]{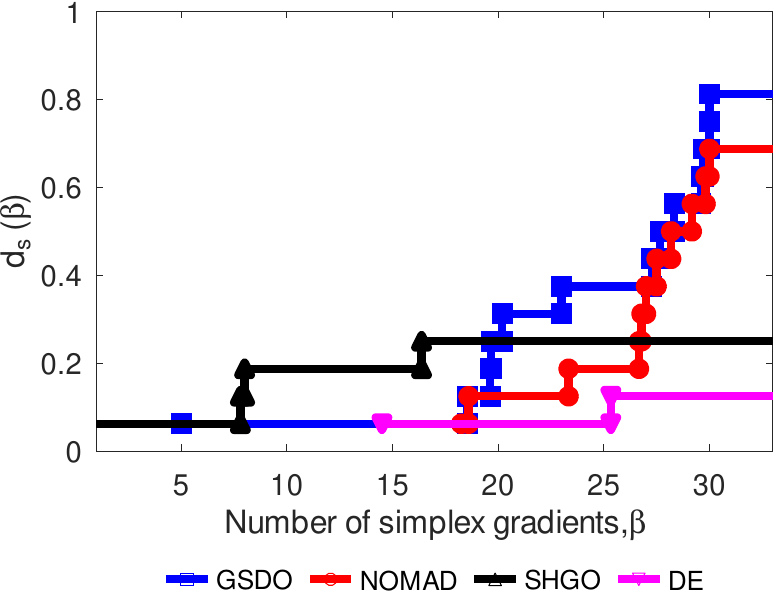}}
	\subfloat[]{ \includegraphics[width=0.4\textwidth]{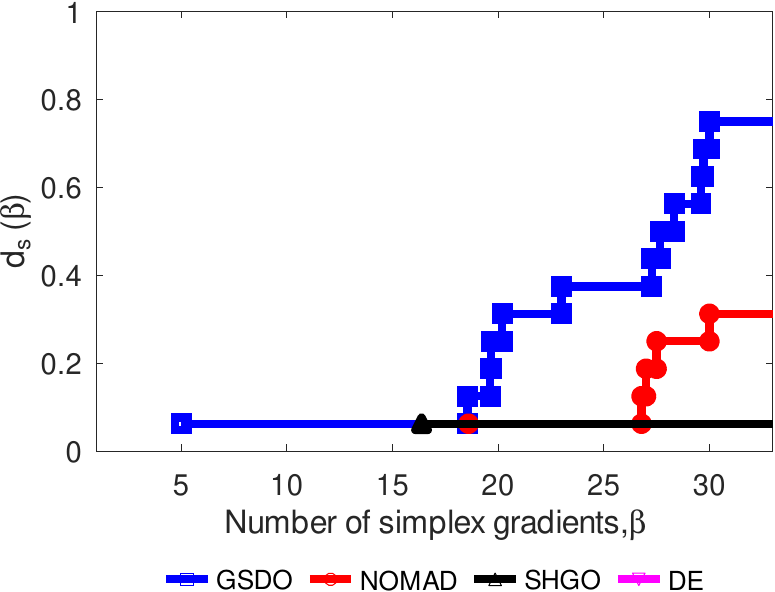}}
	\newline
	\subfloat[]{ \includegraphics[width=0.4\textwidth]{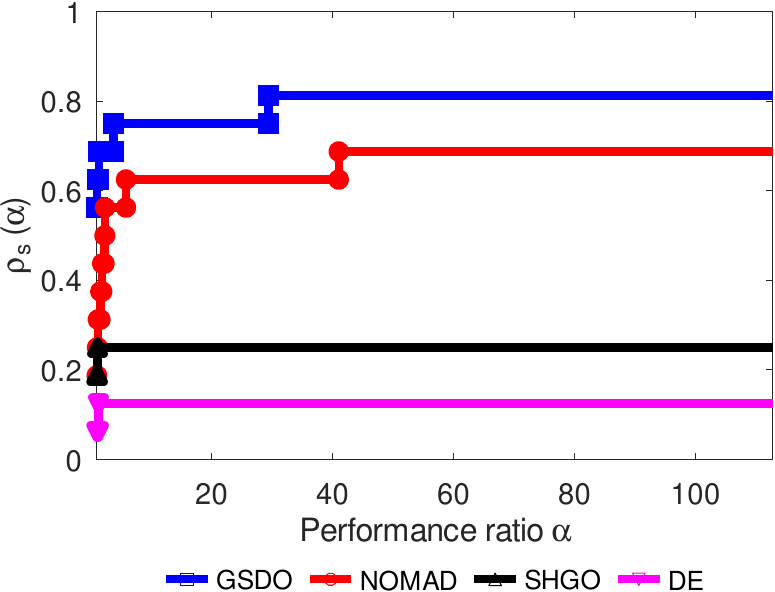}}
	\subfloat[]{ \includegraphics[width=0.4\textwidth]{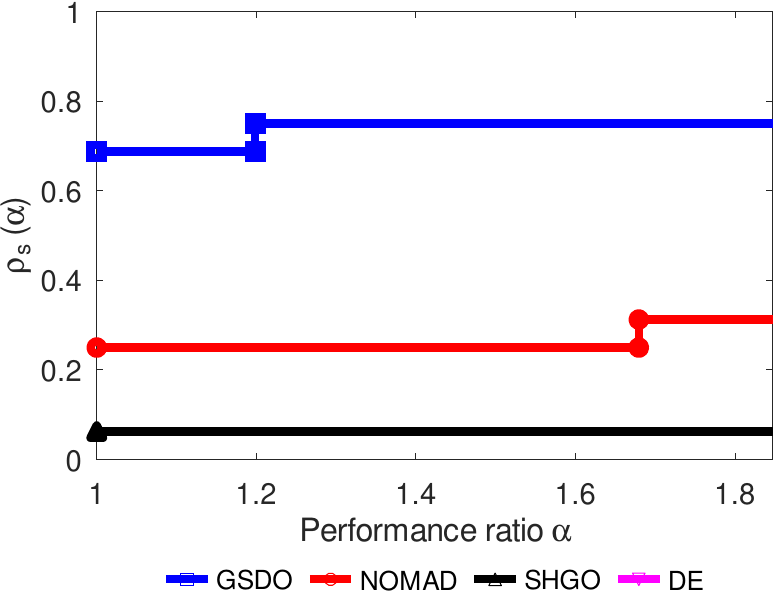}}
	\caption{Set-3: (a) Data profiles at $\tau=0.1$ (b) Data profiles at $\tau=0.01$ (c) Performance profiles at $\tau=0.1$ (d) Performance profiles at $\tau=0.01$}
	\label{fig:dataperfprofile_qusk}
\end{figure}

\begin{figure}[h]
	\subfloat[]{ \includegraphics[width=0.4\textwidth]{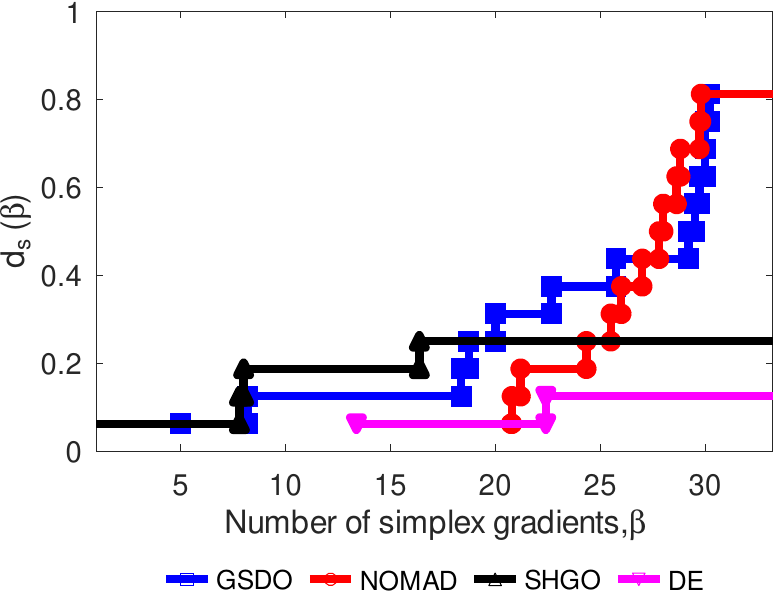}}
	\subfloat[]{ \includegraphics[width=0.4\textwidth]{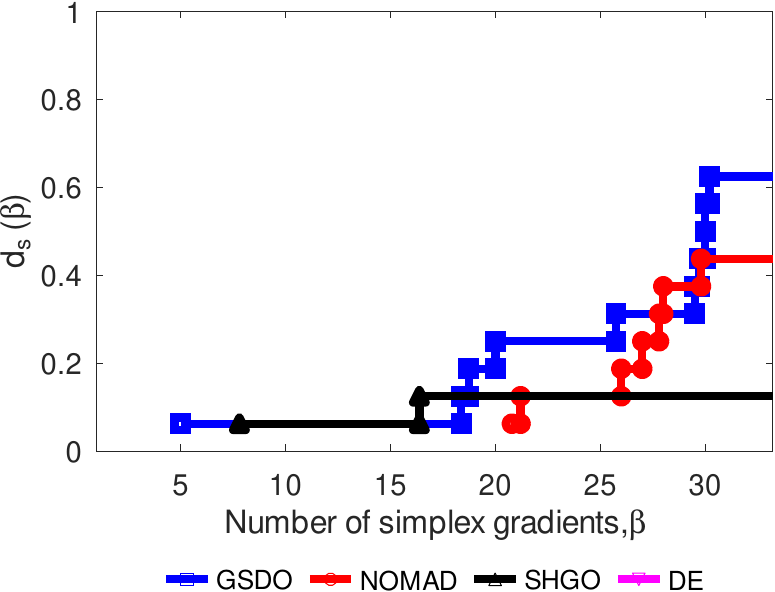}}
	\newline
	\subfloat[]{ \includegraphics[width=0.4\textwidth]{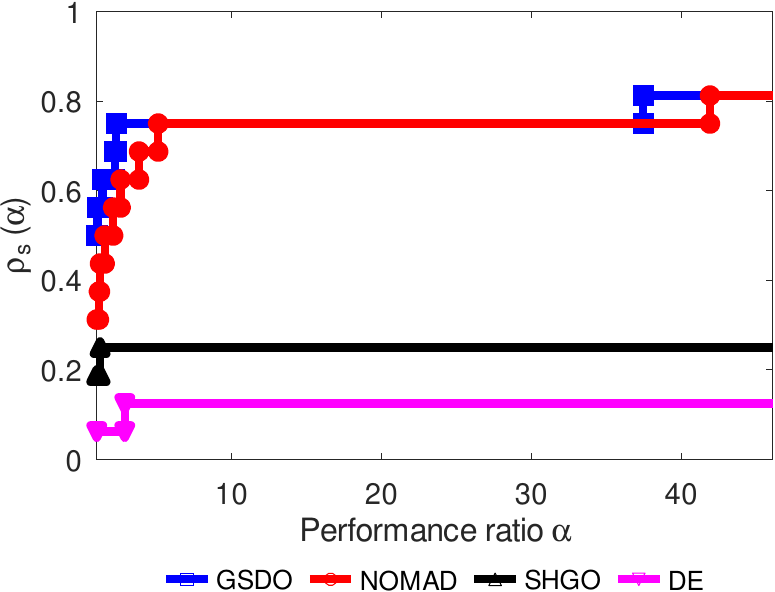}}
	\subfloat[]{ \includegraphics[width=0.4\textwidth]{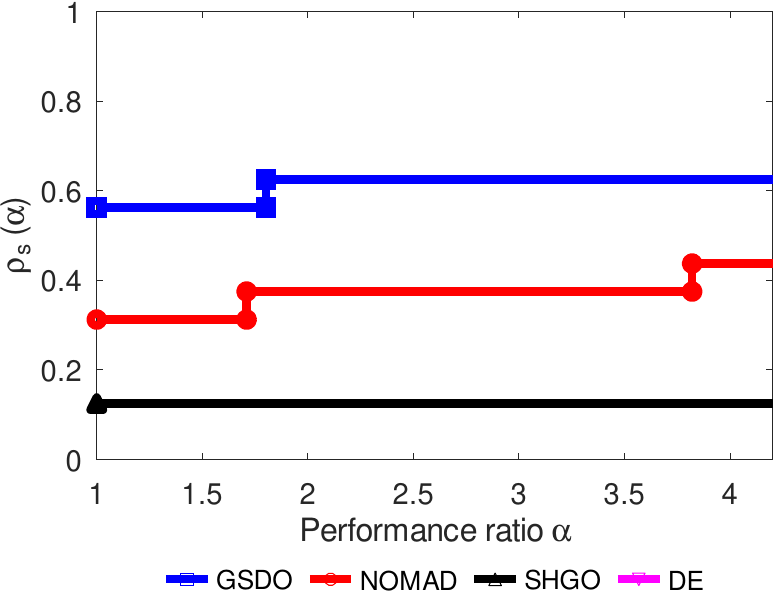}}
	\caption{Set-4: (a) Data profiles at $\tau=0.1$ (b) Data profiles at $\tau=0.01$ (c) Performance profiles at $\tau=0.1$ (d) Performance profiles at $\tau=0.01$}
	\label{fig:dataperfprofile_nusk}
\end{figure}
\clearpage

\begin{table*}[h]
	\center
	\caption {Relative Percentage of Problems Solved}
	\begin{tabular}{|c|c|c|c|c|c|}
		\hline
		&Solver&\multicolumn{2}{c|}{Best Performance} &\multicolumn{2}{c|}{Final Performance}\\
		\hline
		 & &$\tau=0.1$&$\tau=0.01$&$\tau=0.1$&$\tau=0.01$\\
		\hline
		\multirow{4}{*}{Set-1}&GSDO&78&94&94&94\\
		&NOMAD&11&6&28&6\\
		&SHGO&11&0&11&0\\
		&DE&0&0&0&0\\
		\hline
		\multirow{4}{*}{Set-2}&GSDO&62&75&87&75\\
		&NOMAD&12&19&62&31\\
		&SHGO&19&0&25&6\\
		&DE&6&6&19&12\\
		\hline
		\multirow{4}{*}{Set-3}&GSDO&56&69&81&75\\
		&NOMAD&19&25&69&31\\
		&SHGO&19&6&25&6\\
		&DE&6&0&12&0\\
		\hline
		\multirow{4}{*}{Set-4}&GSDO&50&56&81&62\\
		&NOMAD&31&31&81&44\\
		&SHGO&19&12&25&12\\
		&DE&6&0&12&0\\
		\hline
	\end{tabular}
	\label{table:data_perf}
\end{table*}

Next, for better quantitative assessment, we computed relative error with respect to global optimum as defined by \cite{muller2017gosac}. Mathematically, we have
\begin{equation}
	r_e = \left|\frac{f_S - f^*}{f^*}\right|
\end{equation}
where $f_S$ is the optimal feasible solution found by a given solver and $f^*$ is the global optimum (or best known solution).
Figure \ref{fig:errorplot} shows the percentage of problems solved by a specific solver within a relative error $r_e$ ranging from $1\%$ to $10\%$.
Here, for each set, we considered 30 trials for all the test problems in that set, as separate instances.
Hence the total number of problems available were 540 for Set-1 and 480 for the remaining sets.
For relative error less than or equal to $1\%$, NOMAD was able to solve approximately $7\%$ of problems for Set-1, $15\%$ for Set-2, $17\%$ for Set-3 and Set-4.
SHGO and DE failed to solve any problem within $1\%$ relative error.
GSDO, on the other hand, was able to solve approximately $31\%$ for Set-1, $27\%$ for Set-2, $41\%$ for Set-3 and $27\%$ for Set-4.
For $r_e \leq 10\%$, the corresponding numbers for NOMAD, SHGO and DE increased to approximately $17\%$, $11\%$ and $8\%$ for Set-1; $28\%$, $12\%$ and $11\%$ for Set-2; $30\%$, $12\%$ and $11\%$ for Set-3; and $26\%$, $12\%$ and $9\%$ for Set-4.
Similarly, for GSDO, these numbers increased to approximately $45\%$ for Set-1, $35\%$ for Set-2, $44\%$ for Set-3, and $34\%$ for Set-4.
These results show that GSDO consistently performs well in terms of solution quality with respect to actual global optimum.

Another popular metric in DFO to compare performance of different solvers is the evaluation of data and performance profiles, originally proposed by \cite{more2009benchmarking}.
Accordingly, we constructed the performance and data profiles for the test problems using tolerance $\tau=\{0.1,0.01\}$.
Let $\mathcal{P}$ and $\mathcal{S}$ denote the number of test problems and solvers considered.
We denote $w_{s,p}$ as the function evaluations required by the solver $s$ to solve the problem $p$.
The performance ratio $\nu_{s,p}$ is defined as:
\begin{equation}
	\nu_{s,p}=\frac{w_{s,p}}{\text{min}\{w_{s,p}:s\in \mathcal{S}\}}.
\end{equation}
The performance profile for the solver $s$ within performance ratio $\alpha$ is defined as:
\begin{equation}
	\rho_s(\alpha)=\frac{1}{|\mathcal{P}|}\text{size}\{p \in \mathcal{P}:\nu_{s,p} \leq \alpha\}.
\end{equation}
Similarly, the data profile for the solver $s$ for at most $\beta$ simplex gradients is defined as:
\begin{equation}
	d_s(\beta)=\frac{1}{|\mathcal{P}|}\text{size}\{p \in \mathcal{P}:\frac{w_{s,p}}{n_p+1} \leq \beta\},
\end{equation}
where $n_p$ is the number of variables in $p$.

For constructing these profiles, we considered the median of 30 trials as shown in Tables \ref{table:median_qrsk} to \ref{table:median_nusk}.
Parts (a) and (b) of Figures \ref{fig:dataperfprofile_qrsk} to \ref{fig:dataperfprofile_nusk} show the data profiles for all the four solvers at $\tau=0.1$ and $\tau=0.01$.
Parts (c) and (d) of Figure \ref{fig:dataperfprofile_qrsk} to \ref{fig:dataperfprofile_nusk} show the performance profiles (using logarithmic scaling on X-axis) for the four solvers at same tolerance as data profiles.
Table \ref{table:data_perf} reports the performance of each solver relative to the best solver.
The "Best Performance" column refers to $\rho_s(1)$ which indicates the relative percentage of problems on which the solver $s$ performed best within the given accuracy.
The "Final Performance" column refers to the relative percentage of problems solved within the given budget and accuracy for each solver.

From these results, we can infer that the overall performance of GSDO is comparatively better than that of other solvers for all the sets under study.
We now study the performance of these solvers on an application problem.

\subsection{Application Problem: Optimization of Styrene Process}
\begin{table*}[t!]
	\center
	\caption {Results for Styrene Problem}
	\begin{tabular}{|c|c|c|c|c|c|c|c|c|}
		\hline
		&\multicolumn{8}{c|}{Budget}\\
		\hline
		&\multicolumn{4}{c|}{135}&\multicolumn{4}{c|}{270}\\
		\hline
		Solvers&Minimum&Maximum&Median&$N_s$&Minimum&Maximum&Median&$N_s$\\
		\hline
		GSDO&-2.980e+7&-2.144e+7&-2.577e+7&8&-3.138e+7&-1.082e+7&-2.296e+7&18\\
		NOMAD&-3.034e+7&-2.098e+7&-2.977e+7&3&-3.264e+7&-2.333e+7&-2.846e+7&7\\
		SHGO&-&-&-&0&-&-&-&0\\
		DE&-2.636e+7&-7.033e+6&-1.922e+7&12&-3.097e+7&-7.800e+6&-1.962e+7&23\\

		\hline
	\end{tabular}
	\label{table:styrene_results}
\end{table*}
A simulation based model for production of Styrene was developed using the Sequential Modular Simulation (SMS). 
This approach (\cite{peters2003plant}) involves solving of differential equations and other chemical engineering problems in the background.
The optimization problem requires maximization of Net Present Value (NPV) for the process, under the constraints related to meeting industrial and environmental regulations.
A C++ based model and its solution was studied in \cite{AuBeLe08} for this problem.
The model involves 8 variables, all ranging between 0 to 100.
It has 11 constraints, of which 4 are NUSK type and remaining 7 are of QRSK type constraints.
The current best known solution has the objective value of -3.37137e+07.

We experimented at two different budgets, 135 and 270 evaluations.
We performed 30 trials with different starting points for each solver.
The performance of different solvers under consideration, are reported in Table \ref{table:styrene_results}.
A trial is successful for the solver if it is able to find a feasible solution within the given budget.
This is shown under the column ``$N_s$" for the two budgets considered.
The minimum, maximum and median of objective values over these successful trials are also reported in the table for the two budgets.
SHGO was unable to solve the problem for both the budgets.
NOMAD was able to attain the best solution among all solvers, but was successful on only 3 and 7 trials out of 30 for the corresponding budgets.
DE was successful in 12 and 23 trials for the two budgets, but with larger objective values compared to other solvers.
GSDO was successful on 8 and 18 trials, with solution values lying between those of DE and NOMAD.
Thus, performance of GSDO is competitive with other solvers in terms of both solution quality and number of successful trials.

\section{Software Information}\label{software}
The source code for the proposed GSDO solver is available at https://gitlab.com/gcmouli1/gsdo. 
The software is released under GNU General Public License v3.0.
The repository comprises of source files, Makefile for compilation, collection of test problems, and a README file for ready reference.
All instructions regarding the software are provided in the README file.

\section{Conclusion} \label{conclusion}
In this work, a heuristic based open source solver was proposed to address the problem of finding a global solution for the DFO problems with expensive objective and constraint functions.
The proposed approach, GSDO, uses a cubic radial basis function for the surrogate modeling of the given objective and constraint functions.
GSDO is compared with three solvers namely NOMAD, SHGO and DE over 18 test problems and an engineering application problem (Styrene Process Optimization).
These test problems were further investigated in the presence of QRSK, NRSK, QUSK and NUSK types of constraints.
The solvers were first compared on the basis of two parameters, namely solution quality and success percentage out of 30 trials for each problem.
The solution quality was computed in terms of relative error with respect to the best known solutions for the problems considered.
From the results, we infer that GSDO is more successful than other solvers in terms of both these parameters.
Finally, the data and performance profiles were computed, which clearly demonstrate the competitive performance of the proposed method for the given test problems.

\end{document}